\documentclass[11pt,fleqn]{article}
\usepackage{cite}
\usepackage{amsmath, amsthm, amssymb}
\usepackage{bm}
\usepackage{float,afterpage}
\usepackage{leftidx}
\usepackage{srcltx}
\usepackage{graphicx}
\begin{document}

\title{\bf On Selberg-Type Square Matrices Integrals}

\vspace{3.5cm}

\author{{ M. Arashi}\vspace{.5cm}
\\\it
Department of Statistics, School of Mathematical Sciences,\\\it
Shahrood University, Shahrood, Iran\\\it\vspace{-.5cm}\\{\em
Email: m\_arashi\_stat@yahoo.com}}

\date{}
\maketitle

\begin{quotation}
\noindent {\it Abstract:} In this paper we consider Selberg-type
square matrices integrals with focus on Kummer-beta types I \& II
integrals. For generality of the results for real normed division
algebras, the generalized matrix variate Kummer-beta types I \& II
are defined under the abstract algebra. Then Selberg-type
integrals are calculated under orthogonal transformations.
\par

\vspace{9pt} \noindent {\it Key words and phrases:} Selberg-Type
integrals; Real normed division algebras; Kummer-beta
distribution; Random matrix.
\par

\vspace{9pt} \noindent {\it AMS Classification:} Primary: 15B51;
Secondary: 60E05, 05C25.

\end{quotation}\par

\def\blackbox{\ \rule{0.5em}{0.5em}}
\def\tr{\mathop{\rm tr}\nolimits}
\def\Re{\mathop{\rm Re}\nolimits}
\def\Retr{\mathop{\rm Retr}\nolimits}
\def\diag{\mathop{\rm diag}\nolimits}
\def\tr{\mathop{\rm tr}\nolimits}
\def\rank{\mathop{\rm rank}\nolimits}
\def\vec{\mathop{\rm vec}\nolimits}
\def\vecp{\mathop{\rm vecp}\nolimits}
\def\vol{\mathop{\rm Vol}\nolimits}
\def\etr{\mathop{\rm etr}\nolimits}
\def\Rel{\mathop{\rm Re}\nolimits}
\def\var{\mathop{\rm Var}\nolimits}
\def\cov{\mathop{\rm Cov}\nolimits}
\def\corr{\mathop{\rm corr}\nolimits}
\def\func #1{\mathop{\rm #1}\nolimits}%

\newcommand{\balpha}{\boldsymbol{\alpha}}
\newcommand{\bbeta}{\boldsymbol{\beta}}
\newcommand{\bgamma}{\boldsymbol{\gamma}}
\newcommand{\bdelta}{\boldsymbol{\delta}}
\newcommand{\bepsilon}{\boldsymbol{\epsilon}}
\newcommand{\bvarepsilon}{\boldsymbol{\varepsilon}}
\newcommand{\bzeta}{\boldsymbol{\zeta}}
\newcommand{\bet}{\boldsymbol{\eta}}
\newcommand{\btheta}{\boldsymbol{\theta}}
\newcommand{\biota}{\boldsymbol{\iota}}
\newcommand{\bkappa}{\boldsymbol{\kappa}}
\newcommand{\blambda}{\boldsymbol{\lambda}}
\newcommand{\bmu}{\boldsymbol{\mu}}
\newcommand{\bnu}{\boldsymbol{\nu}}
\newcommand{\bxi}{\boldsymbol{\xi}}
\newcommand{\0}{\boldsymbol{0}}
\newcommand{\1}{\boldsymbol{1}}
\newcommand{\bpi}{\boldsymbol{\pi}}
\newcommand{\bvarpi}{\boldsymbol{\varpi}}
\newcommand{\brho}{\boldsymbol{\rho}}
\newcommand{\bvarrho}{\boldsymbol{\varrho}}
\newcommand{\bsigma}{\boldsymbol{\sigma}}
\newcommand{\bvarsigma}{\boldsymbol{\varsigma}}
\newcommand{\btau}{\boldsymbol{\tau}}
\newcommand{\bupsilon}{\boldsymbol{\upsilon}}
\newcommand{\bphi}{\boldsymbol{\phi}}
\newcommand{\bvarphi}{\boldsymbol{\varphi}}
\newcommand{\bchi}{\boldsymbol{\chi}}
\newcommand{\bpsi}{\boldsymbol{\psi}}
\newcommand{\bomega}{\boldsymbol{\omega}}
\newcommand{\bGamma}{\boldsymbol{\Gamma}}
\newcommand{\bDelta}{\boldsymbol{\Delta}}
\newcommand{\bTheta}{\boldsymbol{\Theta}}
\newcommand{\bLambda}{\boldsymbol{\Lambda}}
\newcommand{\bXi}{\boldsymbol{\Xi}}
\newcommand{\bSigma}{\boldsymbol{\Sigma}}
\newcommand{\bUpsilon}{\boldsymbol{\Upsilon}}
\newcommand{\bPhi}{\boldsymbol{\Phi}}
\newcommand{\bPsi}{\boldsymbol{\Psi}}
\newcommand{\bOmega}{\boldsymbol{\Omega}}
\newcommand{\ba}{\boldsymbol{a}}
\newcommand{\bb}{\boldsymbol{b}}
\newcommand{\bA}{\boldsymbol{A}}
\newcommand{\bc}{\boldsymbol{c}}
\newcommand{\bC}{\boldsymbol{C}}
\newcommand{\bd}{\boldsymbol{d}}
\newcommand{\bD}{\boldsymbol{D}}
\newcommand{\be}{\boldsymbol{e}}
\newcommand{\bE}{\boldsymbol{E}}
\newcommand{\bF}{\boldsymbol{F}}
\newcommand{\bg}{\boldsymbol{g}}
\newcommand{\bG}{\boldsymbol{G}}
\newcommand{\bh}{\boldsymbol{h}}
\newcommand{\bH}{\boldsymbol{H}}
\newcommand{\bi}{\boldsymbol{i}}
\newcommand{\bI}{\boldsymbol{I}}
\newcommand{\bj}{\boldsymbol{j}}
\newcommand{\bJ}{\boldsymbol{J}}
\newcommand{\bk}{\boldsymbol{k}}
\newcommand{\bK}{\boldsymbol{K}}
\newcommand{\bl}{\boldsymbol{l}}
\newcommand{\bL}{\boldsymbol{L}}
\renewcommand{\bm}{\boldsymbol{m}}
\newcommand{\bM}{\boldsymbol{M}}
\newcommand{\bn}{\boldsymbol{n}}
\newcommand{\bN}{\boldsymbol{N}}
\newcommand{\bo}{\boldsymbol{o}}
\newcommand{\bO}{\boldsymbol{O}}
\newcommand{\bp}{\boldsymbol{p}}
\newcommand{\bP}{\boldsymbol{P}}
\newcommand{\bq}{\boldsymbol{q}}
\newcommand{\bQ}{\boldsymbol{Q}}
\newcommand{\br}{\boldsymbol{r}}
\newcommand{\bR}{\boldsymbol{R}}
\newcommand{\bs}{\boldsymbol{s}}
\newcommand{\bS}{\boldsymbol{S}}
\newcommand{\bt}{\boldsymbol{t}}
\newcommand{\bT}{\boldsymbol{T}}
\newcommand{\bu}{\boldsymbol{u}}
\newcommand{\bU}{\boldsymbol{U}}
\newcommand{\bv}{\boldsymbol{v}}
\newcommand{\bV}{\boldsymbol{V}}
\newcommand{\bw}{\boldsymbol{w}}
\newcommand{\bW}{\boldsymbol{W}}
\newcommand{\bx}{\boldsymbol{x}}
\newcommand{\bX}{\boldsymbol{X}}
\newcommand{\by}{\boldsymbol{y}}
\newcommand{\bY}{\boldsymbol{Y}}
\newcommand{\bz}{\boldsymbol{z}}
\newcommand{\bZ}{\boldsymbol{Z}}
\newcommand{\A}{\mathcal{A}}
\newcommand{\B}{\mathcal{B}}
\newcommand{\C}{\mathcal{C}}
\newcommand{\D}{\mathcal{D}}
\newcommand{\E}{\mathcal{E}}
\newcommand{\F}{\mathcal{F}}
\newcommand{\G}{\mathcal{G}}
\renewcommand{\H}{\mathcal{H}}
\newcommand{\I}{\mathcal{I}}
\newcommand{\J}{\mathcal{J}}
\newcommand{\K}{\mathcal{K}}
\renewcommand{\L}{\mathcal{L}}
\newcommand{\M}{\mathcal{M}}
\newcommand{\N}{\mathcal{N}}
\renewcommand{\O}{\mathcal{O}}
\renewcommand{\P}{\mathcal{P}}
\newcommand{\Q}{\mathcal{Q}}
\newcommand{\R}{\mathcal{R}}
\renewcommand{\S}{\mathcal{S}}
\newcommand{\T}{\mathcal{T}}
\newcommand{\U}{\mathcal{U}}
\newcommand{\V}{\mathcal{V}}
\newcommand{\W}{\mathcal{W}}
\newcommand{\X}{\mathcal{X}}
\newcommand{\Y}{\mathcal{Y}}
\newcommand{\Z}{\mathcal{Z}}
\def\blackbox{\ \rule{0.5em}{0.5em}}
\def\tr{\mathop{\rm tr}\nolimits}
\def\Re{\mathop{\rm Re}\nolimits}
\def\Retr{\mathop{\rm Retr}\nolimits}
\def\diag{\mathop{\rm diag}\nolimits}
\def\tr{\mathop{\rm tr}\nolimits}
\def\rank{\mathop{\rm rank}\nolimits}
\def\vec{\mathop{\rm vec}\nolimits}
\def\vecp{\mathop{\rm vecp}\nolimits}
\def\vol{\mathop{\rm Vol}\nolimits}
\def\etr{\mathop{\rm etr}\nolimits}
\def\Rel{\mathop{\rm Re}\nolimits}
\def\var{\mathop{\rm Var}\nolimits}
\def\cov{\mathop{\rm Cov}\nolimits}
\def\corr{\mathop{\rm corr}\nolimits}
\def\func #1{\mathop{\rm #1}\nolimits}%

\newtheorem{theorem}{Theorem}[section]
\newtheorem{remark}{Remark}[section]
\newtheorem{definition}{Definition}[section]
\newtheorem{corollary}{Corollary}[theorem]
\newtheorem{lemma}{Lemma}[section]
\renewcommand{\thetheorem}{\arabic{section}.\arabic{theorem}}
\renewcommand{\thedefinition}{\arabic{section}.\arabic{definition}}
\renewcommand{\thelemma}{\arabic{section}.\arabic{lemma}}
\renewcommand{\thesection}{\arabic{section}}
\renewcommand{\thesubsection}{\arabic{section}.\arabic{subsection}}
\renewcommand{\thecorollary}{\arabic{section}.\arabic{theorem}.\arabic{corollary}}

\renewcommand{\thefootnote}{\fnsymbol{footnote}}
\newcommand{\dint}{\int}
\newcommand{\diiint}{\int\int\int}
\newcommand{\dsum}{\sum}

\section{Introduction and Some Preliminaries}
Selberg-type gamma and beta integrals involving scalar functions
of positive definite symmetric matrices as integrand are
considered by several authors, for examples
see~\cite{AskeyRichards} and~\cite{GuptaKabe2005}. Mathai
\cite{Mathai}, (p. 231, 4.1.2) lists Selberg-type gamma integrals
containing a positive signature symmetric matrix. Recently Gupta
and Kabe ~\cite{GuptaKabe2008} presented some results on the
Selberg-type gamma and beta integrals where the integrand is a
scalar function of squared matrix. They also covered skew
symmetric matrices. In this note, we extend the existing results
in the literature for real normed division algebras to cover real,
complex, quaternion and octonion spaces simultaneously. We mainly
focus on Selberg-type square matrices Kummer-gamma and Kummer-beta
integrals.

The hypercomplex multivariate analysis distribution theory
developed by Kabe \cite{Kabe} in order to consider real, complex,
quaternion and octonion spaces simultaneously. Afterward many
researchers extended the results in many directions. Among them,
the works of Prof. Jose A. Diaz-Garcia and his colleague should be
acknowledged. The reader is referred to \cite{Diaz-Garcia2009a},
\cite{Diaz-Garcia2009b}, and \cite{Diaz-Garcia2011}.

A vector-space is always a finite-dimensional module over the
field of real numbers. An algebra $\F$ is a vector space that is
equipped by a bilinear map $m:\F\times\F\rightarrow\F$ termed
multiplication and a non-zero element $1\in\F$ termed the unit
such that $m(1,a)=m(a,1)=a$. As usual abbreviate $m(a,b)=ab$ as
$ab$. Do not assume $\F$ is associative.\\
An algebra $\F$ is a division algebra if given $a,b\in\F$, then
$ab=0$ implies $a=0$ or $b=0$. Equivalently, $\F$ is a division
algebra if the operation of left and right multiplications by any
non-zero element is invertible. A normed division algebra is an
algebra $\F$ that is also a normed vector space with
$\|ab\|=\|a\|\|b\|$. This implies that $\F$ is a division algebra
and $\|1\|=1$.

There are exactly four normed division algebras (according
to~\cite{Baez}):
\begin{enumerate}
\item[(1)] Real Numbers $(\mathbb{R})$,
\item[(2)] Complex Numbers $(\mathbb{C})$,
\item[(3)] Quaternions $(\mathbb{Q})$,
\item[(4)] Octonions $(\mathbb{O})$.
\end{enumerate}
Moreover, they are the only alternative division algebras, and all
division algebras have a real dimension of $1,\;2,\;4$ or $8$,
which is denoted by $\beta$, see \cite{Baez}, theorems 1, 2 and 3.
The parameter $\alpha=2/\beta$
is used, in other mathematical fields, see \cite{EdelmanRao}.\\
Let $\L^\beta_{p,n}$ be the linear space of all $n\times p$
matrices of rank $n\leq p$ over $\F$ with $m$ distinct positive
singular values, where $\F$ denotes a real finite-dimensional
normed division algebra. Let $\F^{n\times p}$ be the set of all
$n\times p$ matrices over $\F$. The dimension of $\F^{n\times p}$
over $\mathbb{R}$ is ¯$np\beta$. Let $\bA\in\F^{n\times p}$, then
$\bA^* = \bar{\bA}^T$ denotes the usual conjugate transpose.\\
The set of matrices $\bH_1\in\F^{n\times p}$ such that
$\bH_1^*\bH_1=\bI_p$ is a manifold denoted $\V^\beta_{p,n}$, is
termed the Stiefel manifold ($\bH_1$ is also known as
semi-orthogonal $(\beta=1)$, semi-unitary $(\beta=2)$,
semi-symplectic $(\beta=4)$ and semi-exceptional type $(\beta=8)$
matrices, see~\cite{DrayManoque}. The dimension of
$\V^\beta_{p,n}$ over $\mathbb{R}$ is
$\left(np\beta-\frac{1}{2}p(p+1)\beta-p\right)$, is the maximal
compact subgroup $\U^\beta(p)$ of $\L^\beta_{p,p}$ and consists of
all matrices $\bH\in\F^{p\times p}$ such that $\bH^*\bH=\bI_p$.
Therefore, $\U^\beta(p)$ is the real orthogonal group $\O(p)$
$(\beta=1)$, the unitary group $\U(p)$ $(\beta=2)$, compact
symplectic group $\S_p(p)$ $(\beta=4)$, or exceptional type
matrices $\O_o(p)$ $(\beta=8)$, for $\F=\mathbb{R}, \mathbb{C},
\mathbb{Q}$ or $\mathbb{O}$, respectively.\\
Denote by $\C^\beta_p$ the real vector space of all
$\bS\in\F^{p\times p}$ such that $\bS=\bS^*$ . Let $\B^\beta_p$ be
the cone of positive definite matrices $\bS\in\F^{p\times p}$;
then $\B^\beta_p$ is an open subset of $\C^\beta_p$. Over
$\mathbb{R}$, $\C^\beta_p$ consists of symmetric matrices; over
$\mathbb{C}$, Hermitian matrices; over $\mathbb{Q}$, quaternionic
Hermitian matrices (also termed self-dual matrices) and over
$\mathbb{O}$, octonionic Hermitian matrices. Generically, the
elements of $\C^\beta_p$ are termed Hermitian matrices,
irrespective of the nature of $\F$. The dimension of $\C^\beta_p$
over $\mathbb{R}$ is $\left(\frac{1}{2}[p(p+1)\beta+p]\right)$.
Let $\D^\beta_p$ be the diagonal subgroup of $\L^\beta_{p,p}$
consisting of all $\bD\in\F^{p\times p}$,
$\bD=\diag(d_1,\ldots,d_p)$.\\
The surface area or volume of the Stiefel manifold
$\V^\beta_{p,n}$ is given by
\begin{eqnarray*}
\mbox{Vol}(\V^\beta_{p,n})=\int_{\bH_1\in\V^\beta_{p,n}}(\bH_1^*d\bH_1)=\frac{2^p\pi^{\frac{np\beta}{2}}}{\Gamma_p^\beta\left(\frac{n\beta}{2}\right)},
\end{eqnarray*}
and therefore
\begin{eqnarray*}
(d\bH_1)=\frac{1}{\mbox{Vol}(\V^\beta_{p,n})}(\bH_1^*d\bH_1)=\frac{\Gamma_p^\beta\left(\frac{n\beta}{2}\right)}{2^p\pi^{\frac{np\beta}{2}}}
(\bH_1^*d\bH_1)
\end{eqnarray*}
is the normalized invariant measure on $\V^\beta_{p,n}$ and
$(d\bH)$, i.e. with $(n=p)$, it defines the normalized Haar
measure on $\U^\beta(p)$ and $\Gamma_p^\beta(a)$ denotes the
multivariate gamma function for the space $\C^\beta_p$, defined by
\begin{eqnarray*}
\Gamma_p^\beta(a)&=&\int_{\bA\in\B_p^\beta}
|\bA|^{a-\frac{(p-1)\beta}{2}-1}\etr(-\bA)d\bA\\
&=&\pi^{\frac{p(p-1)\beta}{4}}\prod_{i=1}^p\Gamma\left(a-\frac{(i-1)\beta}{2}\right),
\end{eqnarray*}
where $\Re(a)>\frac{(p-1)\beta}{2}$, see~\cite{GrossRichards}.

A generalized form of multivariate gamma function is a function of
weight $\kappa$ for the space $\C^\beta_p$ with
$\kappa=(k_1,k_2,\ldots,k_p)$, $k_1\geq k_2\geq\ldots\geq
k_p\geq0$, $\Re(a)\geq(p-1)\beta/2-k_p$, which is defined by:
(see~\cite{GrossRichards} and~\cite{Diaz-Garcia2009a})
\begin{eqnarray*}
\Gamma_p^\beta(a,\kappa)=\int_{\bA\in\B_p^\beta}\etr(-\bA)|\bA|^{a-\frac{(p-1)\beta}{2}-1}q_\kappa(\bA)d\bA=(a)_\kappa^\beta\Gamma_p^\beta(a),
\end{eqnarray*}
where for $\bA\in\C^\beta_p$
\begin{eqnarray*}
q_\kappa(\bA)=|\bA_p|^{k_p}\prod_{i=1}^{p-1}|\bA_i|^{k_i-k_{i+1}}
\end{eqnarray*}
with $\bA_m=(a_{rs})$, $r,s=1,2,\ldots,m$, $m=1,2,\ldots,p$ is
termed the highest weight vector. Also $(a)_\kappa^\beta$ denotes
the generalized Pochhammer symbol of weight $\kappa$, defined by
\begin{eqnarray*}
(a)_\kappa^\beta=\prod_{i=1}^p\left(a-\frac{(i-1)\beta}{2}\right)_{k_i}&=&\frac{\pi^{\frac{p(p-1)\beta}{4}}\prod_{i=1}^p\Gamma\left(a+k_i-\frac{(i-1)\beta}{2}\right)}{\Gamma_p^\beta(a)}\\
&=&\frac{\Gamma_p^\beta(a,\kappa)}{\Gamma_p^\beta(a)},
\end{eqnarray*}
where $\Re(a)>(p-1)\beta/2-k_p$ and
\begin{eqnarray*}
(a)_i=a(a+1)\ldots(a+i-1),
\end{eqnarray*}
is the standard Pochhammer symbol. Thus $\Gamma_p^\beta(a,(0,0,\ldots,0))=\Gamma_p^\beta(a)$\\
The multivariate beta function for the space $\C_p^\beta$ is
defined as (see~\cite{Diaz-Garcia2009a})
\begin{eqnarray*}
B_p^\beta(a,b)&=&\int_{\0<\bX<\bI_p}|\bX|^{a-\frac{(p-1)\beta}{2}-1}|\bI_p-\bX|^{b-\frac{(m-1)\beta}{2}-1}d\bX\\
&=&\int_{\bY\in\B_p^\beta}|\bY|^{a-\frac{(p-1)\beta}{2}-1}|\bI_p+\bY|^{-(a+b)}d\bY\\
&=&\frac{\Gamma_p^\beta(a)\Gamma_p^\beta(b)}{\Gamma_p^\beta(a+b)},
\end{eqnarray*}
where $\bY=(\bI_p-\bX)^{-1}-\bI_p$, $\Re(a),\Re(b)>(p-1)\beta/2$.
From~\cite{KoevEdelman},
\begin{eqnarray*}
\left(\tr(\bX)\right)^k=\sum_\kappa C_\kappa^\beta(\bX),
\end{eqnarray*}
where $C_\kappa^\beta(\bX)$ denotes the zonal polynomials.

Fix complex numbers $a_1,\ldots,a_r$ and $b_1,\ldots,b_s$, and for
all $1\leq i\leq q$ and $1\leq j\leq p$ do not allow
$-b_i+(j-1)\beta/2$ to be a non-negative integer. Then the
hypergeometric function with one matrix argument $_rF_s^\beta$ is
defined to be the real-analytic function on $\C_p^\beta$ given by
the series
\begin{eqnarray*}
_rF_s^\beta(a_1,\ldots,a_r;b_1,\ldots,b_s;\bX)=\sum_{k=0}^\infty\sum_\kappa\frac{(a_1)_\kappa^\beta\ldots(a_r)_\kappa^\beta}
{(b_1)_\kappa^\beta\ldots(b_s)_\kappa^\beta}\;\frac{C_\kappa^\beta(\bX)}{k!}.
\end{eqnarray*}
For convergence properties see~\cite{GrossRichards}.\\
Then we have
\begin{eqnarray*}
_0F_0^\beta(\bX)=\sum_{k=0}^\infty \sum_\kappa
\frac{C_\kappa^\beta(\bX)}{k!}=\etr(\bX).
\end{eqnarray*}
For more details and results the interested reader is referred
to~\cite{Diaz-Garcia2011},~\cite{Dimitriu},
and~\cite{GrossRichards}.

Analogous to ~\cite{Muirhead}, define the confluent hypergeometric
function $\Psi^\beta$ of $p\times p$ matrix $\bX\in\B_p^\beta$ as
\begin{equation*}
\Psi^\beta(a,c;\bX)=\frac{1}{\Gamma_p^\beta(a)}\int_{\bY\in\B_p^\beta}\etr(-\bX\bY)|\bY|^{a-\frac{(p-1)\beta}{2}-1}
|\bI_p+\bY|^{c-a-\frac{(p-1)\beta}{2}-1}\textnormal{d}\bY,
\end{equation*}
where $\Re(a)>(p-1)\beta/2$.

\section{Matrix Variate Distributions}
In this study we consider some matrix variate distribution for our
purpose. Consider the following definitions from
\cite{Diaz-Garcia2009a} and \cite{Diaz-Garcia2009b}.
\begin{definition} Let $\bX\in\L_{p,n}^\beta$ be a random
matrix, and $\bSigma\in\B_p^\beta$ and $\bTheta\in\B_n^\beta$ be
parameter matrices.
\begin{enumerate}
\item[1.] (\textbf{Matrix Variate Normal} ) The random matrix $\bX$ is said to have a matrix variate normal
distribution denoted by $\bX\sim N_{n\times
p}^\beta(\bmu,\bSigma,\bTheta)$, with mean $\bmu$ and
$\textnormal{Cov}(\vec \bX^*)=\bTheta\otimes\bSigma$, if its
density function is given by
\begin{eqnarray*}
\frac{\beta^{\frac{\beta pn}{2}}}{(2\pi)^{\frac{\beta
pn}{2}}|\bSigma|^{\frac{\beta n}{2}}|\bTheta|^{\frac{\beta
p}{2}}}\;\etr\left\{-\frac{\beta}{2}\;\bSigma^{-1}(\bX-\bmu)^*\bTheta^{-1}(\bX-\bmu)\right\}.
\end{eqnarray*}
\item[2.] (\textbf{Wishart}) Let $\bX\sim N_{n\times p}^\beta(\0,\bSigma,\bI_n)$ and define
$\bS=\bX^*\bX$, then $\bS$ is said to have a central Wishart
distribution $\bS\sim W_p^\beta(n,\bSigma)$ with $n$ degrees of
freedom and parameter matrix $\bSigma$. Moreover, its density
function is given by
\begin{eqnarray*}
\frac{\beta^{\frac{\beta pn}{2}}}{2^{\frac{\beta
pn}{2}}\Gamma_p^\beta\left(\frac{\beta
n}{2}\right)|\bSigma|^{\frac{\beta
n}{2}}}\;|\bS|^{\frac{\beta(n-p+1)}{2}-1}\etr\left\{-\frac{\beta}{2}\;\bSigma^{-1}\bS\right\},
\end{eqnarray*}
where $n\geq p-1$.
\item[3.] (\textbf{Matrix Variate T Type II}) The random matrix $\bX$ is said to have a matrix variate $T$-distribution of type II
denoted by $\bX\sim T_{n,p}^\beta(\nu,\bI_p)$, if its density
function is given by
\begin{eqnarray*}
\frac{\Gamma_p^\beta\left(\frac{\beta(n+\nu)}{2}\right)}{\pi^{\frac{\beta
pn}{2}}\Gamma_p^\beta\left(\frac{\beta\nu}{2}\right)}\;|\bI_p+\bX^*\bX|^{-\frac{\beta(n+\nu)}{2}},\quad
\nu>p.
\end{eqnarray*}
This distribution is also termed as matrix variate Pearson Type
VII distribution.
\item[4.] (\textbf{Gegenbauer Type II}) The random matrix $\bX$ is said to have a Gegenbauer Type
II distribution $\bX\sim G_{n,p}^\beta(\nu,\bI_p)$, if its density
function is given by
\begin{eqnarray*}
\frac{\Gamma_p^\beta\left(\frac{\beta(n+\nu)}{2}\right)}{\pi^{\frac{\beta
pn}{2}}\Gamma_p^\beta\left(\frac{\beta\nu}{2}\right)}\;|\bI_p-\bX^*\bX|^{\frac{\beta(\nu-p+1)}{2}-1},\quad
\nu>p-1.
\end{eqnarray*}
This distribution is known in statistical bibliography as the
matrix variate inverted T or matrix variate Pearson Type II
distribution.
\item[5.] (\textbf{T-Laguerre Type II Ensemble}) The random matrix $\bX$ is said to have
a T-Laguerre Type II ensemble distribution $\bX\sim
TL_{n,p}^\beta(\nu)$, if its density function is given by
\begin{eqnarray*}
\frac{1}{B_p^\beta\left(\frac{\beta \nu}{2},\frac{\beta
n}{2}\right)}\;|\bX^*\bX|^{\frac{\beta(n-p+1)}{2}-1}|\bI_p+\bX^*\bX|^{-\frac{\beta(n+\nu)}{2}},
\end{eqnarray*}
where $\nu\geq p-1,\quad n\geq p-1$.

This distribution is also known as the Studentized Wishart
distribution
\item[6.] (\textbf{Gegenbauer-Laguerre Type II Ensemble}) The random matrix $\bX$ is said to have
a Gegenbauer-Laguerre type II ensemble distribution $\bX\sim
GL_{n,p}^\beta(\nu)$, if its density function is given by
\begin{eqnarray*}
\frac{1}{B_p^\beta\left(\frac{\beta \nu}{2},\frac{\beta
n}{2}\right)}\;|\bX^*\bX|^{\frac{\beta(n-p+1)}{2}-1}|\bI_p-\bX^*\bX|^{\frac{\beta(\nu-p+1)}{2}-1},
\end{eqnarray*}
where $\nu\geq p-1,\quad n\geq p-1$.
\end{enumerate}
\end{definition}
The following result directly obtains from Definition 2.1 and
Proposition 4 of~\cite{Diaz-Garcia2011}.
\begin{theorem}~
\begin{enumerate}
\item[(1)] Let $\bX\sim T_{n,p}^\beta(\nu,\bI_p)$. Define
$\bS_1=\bX^*\bX\in\B_p^\beta$. Then $\bS_1$ has the following
density
\begin{eqnarray*}
\frac{1}{B_p^\beta\left(\frac{\beta \nu}{2},\frac{\beta
n}{2}\right)}\;|\bS_1|^{\frac{\beta(n-p+1)}{2}-1}|\bI_p+\bS_1|^{-\frac{\beta(n+\nu)}{2}},
\end{eqnarray*}
where $n\geq p$.
\item[(2)] Let $\bY\sim G_{n,p}^\beta(\nu,\bI_p)$. Define
$\bS_2=\bY^*\bY\in\B_p^\beta$. Then $\bS_2$ has the following
density
\begin{eqnarray*}
\frac{1}{B_p^\beta\left(\frac{\beta \nu}{2},\frac{\beta
n}{2}\right)}\;|\bS_2|^{\frac{\beta(n-p+1)}{2}-1}|\bI_p-\bS_2|^{\frac{\beta(\nu-p+1)}{2}-1},
\end{eqnarray*}
where $n\geq p$.
\end{enumerate}
\end{theorem}
The following definitions are analogous to the results
of~\cite{Nagar2002} and~\cite{Nagar2001} respectively for real
normed division algebras.
\begin{definition} (\textbf{Kummer-Beta Type I}) The random matrix $\bX$ is said to have
a Kummer-beta type I distribution $\bX\sim
KB1_p^\beta(\alpha_1,\alpha_2,\bSigma)$, where
$\bSigma\in\B_p^\beta$ if its density function is given by
\begin{eqnarray*}
K_1(\alpha_1,\alpha_2,\beta,\bSigma)\etr\left(-\bSigma\bX\right)|\bX|^{\alpha_1-\frac{(p-1)\beta}{2}-1}|\bI_p-\bX|^{\alpha_2-\frac{(p-1)\beta}{2}-1},
\end{eqnarray*}
where $\0<\bX<\bI_p$, $\alpha_1\geq(p-1)\beta/2$ and
$\alpha_2\geq(p-1)\beta/2$. The normalizing constant is given by
\begin{eqnarray*}
\left\{K_1(\alpha_1,\alpha_2,\beta,\bSigma)\right\}^{-1}&&\\
=\int_{\0<\bY<\bI_p}\etr\left(-\bSigma\bY\right)|\bY|^{\alpha_1-\frac{(p-1)\beta}{2}-1}|\bI_p-\bY|^{\alpha_2-\frac{(p-1)\beta}{2}-1}\textnormal{d}\bY&&\\
=B_p^\beta\left(\alpha_1,\alpha_2\right)\;_1F_1^\beta(\alpha_1,\alpha_1+\alpha_2;-\bSigma).&&
\end{eqnarray*}
\end{definition}
\begin{definition} (\textbf{Kummer-Beta Type II}) The random matrix $\bX$ is said to have
a Kummer-beta type II distribution $\bX\sim
KB2_p^{\beta}(\alpha_1,\alpha_2,\bSigma)$, where
$\bSigma\in\B_p^\beta$ if its density function is given by
\begin{eqnarray*}
K_2(\alpha_1,\alpha_2,\beta,\bSigma)\etr\left(-\bSigma\bX\right)|\bX|^{\alpha_1-\frac{(p-1)\beta}{2}-1}|\bI_p+\bX|^{-\alpha_2},\quad\bX>\0,
\end{eqnarray*}
where using the confluent hypergeometric function, the normalizing
constant is given by
\begin{eqnarray*}
\left\{K_2(\alpha_1,\alpha_2,\beta,\bSigma)\right\}^{-1}
&=&\int_{\bY\in\B_p^\beta}\etr\left(-\bSigma\bY\right)|\bY|^{\alpha_1-\frac{(p-1)\beta}{2}-1}|\bI_p+\bY|^{-\alpha_2}\textnormal{d}\bY\\
&=&\Gamma_p^\beta(\alpha_1)\Psi^\beta\left(\alpha_1,\alpha_1-\alpha_2\frac{(p-1)\beta}{2}+1;\bSigma\right).
\end{eqnarray*}
\end{definition}
The following result is an extension to Theorem 3.1
of~\cite{Nagar2002} for real normed division algebras.
\begin{lemma}
Let $\bU\sim KB1_p^\beta(\alpha_1,\alpha_2,\bSigma)$. Then for the
given matrices $\bPsi\in\C_p^\beta$, $\bOmega\in\B_p^\beta$ and
$\bOmega-\bPsi\in\B_p^\beta$ $(\bOmega>\bPsi)$, the random matrix
$\bX$ defined by
\begin{eqnarray*}
\bX=(\bOmega-\bPsi)^{\frac{1}{2}}\bU(\bOmega-\bPsi)^{\frac{1}{2}}+\bPsi
\end{eqnarray*}
has generalized matrix variate Kummer-beta type I distribution
with the following density function
\begin{eqnarray*}
C_1(\alpha_1,\alpha_2,\beta,\bTheta,\bOmega,\bPsi)\etr(-\bTheta\bX)|\bX-\bPsi|^{\alpha_1-\frac{(p-1)\beta}{2}-1}|\bOmega-\bX|^{\alpha_2-\frac{(p-1)\beta}{2}-1},
\end{eqnarray*}
where $\bPsi<\bX<\bOmega$,
$\bTheta=(\bOmega-\bPsi)^{-\frac{1}{2}}\bSigma(\bOmega-\bPsi)^{-\frac{1}{2}}$
and
\begin{eqnarray*}
C_1(\alpha_1,\alpha_2,\beta,\bTheta,\bOmega,\bPsi)&=&K_1(\alpha_1,\alpha_2,\beta,(\bOmega-\bPsi)^{\frac{1}{2}}\bTheta(\bOmega-\bPsi)^{\frac{1}{2}})\\
&&\times\etr(\bTheta\bPsi)
|\bOmega-\bPsi|^{-(\alpha_1+\alpha_2)+\frac{(p-1)\beta}{2}+1}.
\end{eqnarray*}
In this case we use the notation $\bX\sim
GKB1_p^\beta(\alpha_1,\alpha_2,\bTheta,\bOmega,\bPsi)$.
\end{lemma}
\noindent\textbf{Proof:} The proof directly follows by Definition
2.2 and the fact that the Jacobian of transformation is
$J(\bU\rightarrow\bX)=|\bOmega-\bPsi|^{-\frac{\beta(p-1)}{2}-1}$.

The following result is an extension to Theorem 2.2
of~\cite{Nagar2001} for real normed division algebras.
\begin{lemma}
Let $\bU\sim KB2_p^\beta(\alpha_1,\alpha_2,\bSigma)$. Then for the
given matrices $\bPsi\in\C_p^\beta$, $\bOmega\in\B_p^\beta$ and
$\bOmega+\bPsi\in\B_p^\beta$, the random matrix $\bX$ defined by
\begin{eqnarray*}
\bX=(\bOmega+\bPsi)^{\frac{1}{2}}\bU(\bOmega+\bPsi)^{\frac{1}{2}}+\bPsi
\end{eqnarray*}
has generalized matrix variate Kummer-beta type II with the
following density function
\begin{eqnarray*}
C_2(\alpha_1,\alpha_2,\beta,\bTheta,\bOmega,\bPsi)\etr(-\bTheta\bX)|\bX-\bPsi|^{\alpha_1-\frac{(p-1)\beta}{2}-1}|\bOmega+\bX|^{-\alpha_2},
\end{eqnarray*}
where $\bPsi<\bX$,
$\bTheta=(\bOmega+\bPsi)^{-\frac{1}{2}}\bSigma(\bOmega+\bPsi)^{-\frac{1}{2}}$
and
\begin{eqnarray*}
C_2(\alpha_1,\alpha_2,\beta,\bTheta,\bOmega,\bPsi)&=&K_2(\alpha_1,\alpha_2,\beta,(\bOmega+\bPsi)^{\frac{1}{2}}\bTheta(\bOmega+\bPsi)^{\frac{1}{2}})\\
&&\times\etr(\bTheta\bPsi)|\bOmega+\bPsi|^{-\alpha_1+\alpha_2}.
\end{eqnarray*}
In this case we use the notation $\bX\sim
GKB2_p^\beta(\alpha_1,\alpha_2,\bTheta,\bOmega,\bPsi)$.
\end{lemma}
\noindent\textbf{Proof:} The proof directly follows by Definition
3.3 and the fact that the Jacobian of transformation is
$J(\bU\rightarrow\bX)=|\bOmega+\bPsi|^{-\frac{\beta(p-1)}{2}-1}$.

\setcounter{equation}{0}
\section{Selberg Type Square Matrices Integrals}
In this section we are interested in evaluating the integrals of
the form
\begin{eqnarray}\label{eq1}
\int h(\bLambda)\prod_{i<j}^p(\lambda_i-\lambda_j)^\beta
\textnormal{d}\bLambda,
\end{eqnarray}
where $\bLambda=\textnormal{Diag}(\lambda_1,\ldots,\lambda_p)$,
$\lambda_1>\ldots>\lambda_p>0$, for some function $h$. In this
regard, we need the following essential result due
to~\cite{Diaz-Garcia2009a}. Note that similar results can also be
found in~\cite{Diaz-Garcia2009b} for general form of the
distributions defined in Definition 2.1. However the purpose of
this study is to evaluate integrals of the form \eqref{eq1}.
\begin{theorem} Let $\bX\in\B_p^\beta$ be a
random matrix with density function $f(\bX)$. Then the joint
density function of the eigenvalues $\lambda_1,\ldots,\lambda_p$
of $\bX$ is
\begin{eqnarray*}
g(\lambda_1,\ldots,\lambda_p)&=&\frac{\pi^{\frac{1}{2}p^2\beta+\varrho}}{\Gamma_p^\beta\left(\frac{p\beta}{2}\right)}\;\prod_{i<j}^p(\lambda_i-\lambda_j)^\beta
\int_{\bH\in\U^\beta(p)} f(\bH\bLambda\bH^*)(\textnormal{d}\bH),
\end{eqnarray*}
where $(\textnormal{d}\bH)$ is the normalized Haar measure and
\begin{eqnarray*}
\varrho=\left\{\begin{array}{cc}
      0, & \beta=1;\\
      -p, & \beta=2;\\
      -2p, & \beta=4;\\
      -4p, & \beta=8.\end{array}\right.
\end{eqnarray*}
\end{theorem}

\begin{lemma} Let $\bX\in\B_p^\beta$ be a
random matrix with density function $f(\bX)$. Then we have
\begin{eqnarray*}
\int_{\bH\in\U^\beta(p)}\left(\int_{\bLambda}f(\bH\bLambda\bH^*)\prod_{i<j}^p(\lambda_i-\lambda_j)^\beta
\textnormal{d}\bLambda\right)
(\textnormal{d}\bH)=\frac{\Gamma_p^\beta\left(\frac{p\beta}{2}\right)}{\pi^{\frac{1}{2}p^2\beta+\varrho}}.
\end{eqnarray*}
\end{lemma}
\noindent\textbf{Proof}: Since $\int_{\bLambda}
g(\lambda_1,\ldots,\lambda_p)\textnormal{d}\bLambda=1$, by making
use of Theorem 3.1 and changing the order of integration, the
result follows.

In sequel we proceed by giving some examples of gamma and beta
integrals for real normed division algebras.
\subsection{Examples}
In this part, we give some examples of Selberg-type integrals.\\

\textbf{Example 1.} Based on Definition 2.1, and the fact that
$\bH^*\bH=\bI_p$, we have the following results using Lemma 3.1
\begin{enumerate}
\item[(1:1)] Suppose that $\bX\sim W_p^\beta(n,\bI_p)$, then using the density of $\bS$ given in Definition 2.1(2), we have (extension to \cite{Mathai})
\begin{eqnarray*}
\int_{\bLambda}|\bLambda|^{\frac{\beta(n-p+1)}{2}-1}\etr\left\{-\frac{\beta}{2}\;\bLambda\right\}\prod_{i<j}^p(\lambda_i-\lambda_j)^\beta
\textnormal{d}\bLambda=\frac{2^{\frac{\beta
pn}{2}}\left[\Gamma_p^\beta\left(\frac{p\beta}{2}\right)\right]^2}{\beta^{\frac{\beta
pn}{2}}\pi^{\frac{1}{2}p^2\beta+\varrho}}
\end{eqnarray*}
\item[(1:2)] Suppose that $\bX\sim T_{n,p}^\beta(\nu,\bI_p)$, then from the density of $\bS_1=\bX^*\bX$ given in Definition 2.1(3), we have (extension to \cite{Gupta Kabe2005})
\begin{eqnarray*}
\int_{\bLambda}|\bLambda|^{\frac{\beta(n-p+1)}{2}-1}|\bI_p+\bLambda|^{-\frac{\beta(n+\nu)}{2}}
\prod_{i<j}^p(\lambda_i-\lambda_j)^\beta
\textnormal{d}\bLambda=\frac{\Gamma_p^\beta\left(\frac{p\beta}{2}\right)
B_p^\beta\left(\frac{\beta \nu}{2},\frac{\beta
n}{2}\right)}{\pi^{\frac{1}{2}p^2\beta+\varrho}}.
\end{eqnarray*}
In the same fashion
\item[(1:3)] Suppose that $\bX\sim G_{n,p}^\beta(\nu,\bI_p)$, then
from the density of $\bS_2=\bX^*\bX$ given in Definition 2.1(4),
we have
\begin{eqnarray*}
\int_{\bLambda}|\bLambda|^{\frac{\beta(n-p+1)}{2}-1}|\bI_p-\bLambda|^{\frac{\beta(\nu-p+1)}{2}-1}\prod_{i<j}^p(\lambda_i-\lambda_j)^\beta
\textnormal{d}\bLambda=\frac{\Gamma_p^\beta\left(\frac{p\beta}{2}\right)
B_p^\beta\left(\frac{\beta \nu}{2},\frac{\beta
n}{2}\right)}{\pi^{\frac{1}{2}p^2\beta+\varrho}}.
\end{eqnarray*}
\end{enumerate}
\indent\textbf{Example 2.}
\begin{enumerate}
\item[(2:1)] Suppose that $\bX\sim KB1_p^\beta(\alpha_1,\alpha_2,\bSigma)$, then by
Definition 2.2 we get
\begin{eqnarray}\label{eq32}
\int_{\bH}\int_{\bLambda}\etr\left(\bSigma\bH\bLambda\bH^*\right)|\bI_p-\bLambda|^{\alpha_1-\frac{(p-1)\beta}{2}-1}|\bLambda|^{\alpha_2-\frac{(p-1)\beta}{2}-1}
\prod_{i<j}^p(\lambda_i-\lambda_j)^\beta
\textnormal{d}\bLambda(\textnormal{d}\bH)&&\nonumber\\
=\frac{\Gamma_p^\beta\left(\frac{p\beta}{2}\right)B_p^\beta\left(\alpha_1,\alpha_2\right)}
{\pi^{\frac{1}{2}p^2\beta+\varrho}}\;_1F_1^\beta(\alpha_1,\alpha_1+\alpha_2;-\bSigma).&&
\end{eqnarray}
It is easily seen that taking $\bSigma=\0$,
$\alpha_1=\frac{\nu\beta}{2}$ and $\alpha_2=\frac{n\beta}{2}$ in
\eqref{eq32} gives item (1:3) in the above. Further, for the case
$\bSigma=\bI_p$, we have that
\begin{eqnarray*}
\int_{\bLambda}\etr\left(\bLambda\right)|\bI_p-\bLambda|^{\alpha_1-\frac{(p-1)\beta}{2}-1}|\bLambda|^{\alpha_2-\frac{(p-1)\beta}{2}-1}
\prod_{i<j}^p(\lambda_i-\lambda_j)^\beta
\textnormal{d}\bLambda&&\nonumber\\
=\frac{\Gamma_p^\beta\left(\frac{p\beta}{2}\right)B_p^\beta\left(\alpha_1,\alpha_2\right)}
{\pi^{\frac{1}{2}p^2\beta+\varrho}}\;_1F_1^\beta(\alpha_1,\alpha_1+\alpha_2;-\bI_p).&&
\end{eqnarray*}
Note that the integral in \eqref{eq32} is not of required form,
thus making use of (see \cite{Diaz-Garcia2009a})
\begin{eqnarray}
\int_{\bH}\etr\left(\bSigma\bH\bLambda\bH^*\right)(\textnormal{d}\bH)&=&
\sum_{k=0}^\infty\frac{1}{k!}\int_{\bH}C_\kappa(\bSigma\bH\bLambda\bH^*)(\textnormal{d}\bH)\cr
&=&\sum_{k=0}^\infty\frac{C_\kappa(\bSigma)C_\kappa(\bLambda)}{k!C_\kappa(\bI_p)}
\end{eqnarray}
we obtain
\begin{eqnarray}
\sum_{k=0}^\infty\frac{C_\kappa(\bSigma)}{k!C_\kappa(\bI_p)}
\int_{\bLambda}|\bI_p-\bLambda|^{\alpha_1-\frac{(p-1)\beta}{2}-1}|\bLambda|^{\alpha_2-\frac{(p-1)\beta}{2}-1}\nonumber\\
\times C_\kappa(\bLambda)\prod_{i<j}^p(\lambda_i-\lambda_j)^\beta
\textnormal{d}\bLambda&&\nonumber\\
=\frac{\Gamma_p^\beta\left(\frac{p\beta}{2}\right)B_p^\beta\left(\alpha_1,\alpha_2\right)}
{\pi^{\frac{1}{2}p^2\beta+\varrho}}\;_1F_1^\beta(\alpha_1,\alpha_1+\alpha_2;-\bSigma).&&
\end{eqnarray}
\item[(2:2)] Suppose that $\bX\sim KB2_p^\beta(\alpha_1,\alpha_2,\bSigma)$, then by
Definition 2.3 we get
\begin{eqnarray}\label{eq2}
\int_{\bH}\int_{\bLambda}\etr\left(\bSigma\bH\bLambda\bH^*\right)|\bLambda|^{\alpha_1-\frac{(p-1)\beta}{2}-1}|\bI_p+\bLambda|^{-\alpha_2}
\prod_{i<j}^p(\lambda_i-\lambda_j)^\beta
\textnormal{d}\bLambda(\textnormal{d}\bH)&&\nonumber\\
=\frac{\Gamma_p^\beta\left(\frac{p\beta}{2}\right)\Gamma_p^\beta(\alpha_1)}
{\pi^{\frac{1}{2}p^2\beta+\varrho}}\Psi^\beta\left(\alpha_1,\alpha_1-\alpha_2\frac{(p-1)\beta}{2}+1;\bSigma\right).&&
\end{eqnarray}
It is also interesting to see that for the case $\bSigma=\bI_p$ we
get
\begin{eqnarray*}
\int_{\bLambda}\etr\left(\bLambda\right)|\bLambda|^{\alpha_1-\frac{(p-1)\beta}{2}-1}|\bI_p+\bLambda|^{-\alpha_2}
\prod_{i<j}^p(\lambda_i-\lambda_j)^\beta
\textnormal{d}\bLambda&&\nonumber\\
=\frac{\Gamma_p^\beta\left(\frac{p\beta}{2}\right)\Gamma_p^\beta(\alpha_1)}
{\pi^{\frac{1}{2}p^2\beta+\varrho}}\Psi^\beta\left(\alpha_1,\alpha_1-\alpha_2\frac{(p-1)\beta}{2}+1;\bI_p\right).&&
\end{eqnarray*}
\item[(2:3)] Suppose that $\bX\sim
GKB1_p^\beta(\alpha_1,\alpha_2,\bTheta,\bOmega,\bPsi)$, then by
Lemma 2.1 we have
\begin{eqnarray}\label{eq3}
\int_{\bH}\int_{\bLambda}\etr\left(-\bTheta\bH\bLambda\bH^*\right)
|\bH\bLambda\bH
^*-\bPsi|^{\alpha_1-\frac{(p-1)\beta}{2}-1}&&\nonumber\\|\bOmega-\bH\bLambda\bH^*|^{\alpha_2-\frac{(p-1)\beta}{2}-1}
\prod_{i<j}^p(\lambda_i-\lambda_j)^\beta
\textnormal{d}\bLambda(\textnormal{d}\bH)&&\nonumber\\
=\frac{\Gamma_p^\beta\left(\frac{p\beta}{2}\right)B_p^\beta\left(\alpha_1,\alpha_2\right)}
{\pi^{\frac{1}{2}p^2\beta+\varrho}}\;_1F_1^\beta(\alpha_1,\alpha_1+\alpha_2;-\bSigma)\nonumber\\
\times\etr(-\bTheta\bPsi)
|\bOmega-\bPsi|^{\alpha_1+\alpha_2-\frac{(p+1)\beta}{2}-1},&&
\end{eqnarray}
where
$\bSigma=(\bOmega-\bPsi)^{\frac{1}{2}}\bTheta(\bOmega-\bPsi)^{\frac{1}{2}}$.
\item[(2:4)] Suppose that $\bX\sim
GKB2_p^\beta(\alpha_1,\alpha_2,\bTheta,\bOmega,\bPsi)$, then by
Lemma 2.2 we have
\begin{eqnarray*}\label{eq3}
\int_{\bH}\int_{\bLambda}\etr\left(-\bTheta\bH\bLambda\bH^*\right)
|\bH\bLambda\bH^*-\bPsi|^{\alpha_1-\frac{(p-1)\beta}{2}-1}|\bOmega+\bH\bLambda\bH^*|^{-\alpha_2}
\prod_{i<j}^p(\lambda_i-\lambda_j)^\beta
\textnormal{d}\bLambda(\textnormal{d}\bH)&&\cr
=\frac{\Gamma_p^\beta\left(\frac{p\beta}{2}\right)\Gamma_p^\beta\left(\alpha_1\right)}
{\pi^{\frac{1}{2}p^2\beta+\varrho}}\etr(-\bTheta\bPsi)|\bOmega+\bPsi|^{\alpha_1-\alpha_2}&&\cr
\times
\Psi^\beta\left(\alpha_1,\alpha_1-\alpha_2\frac{(p-1)\beta}{2}+1;(\bOmega+\bPsi)^{\frac{1}{2}}\bTheta(\bOmega+\bPsi)^{\frac{1}{2}}\right).
\end{eqnarray*}

\end{enumerate}

\section*{Acknowledgments}
The author would like to thank the anonymous referees for their
constructive comments that let to put many details on the paper
and improved the presentation.

\end{document}